\documentclass[twoside,11pt]{article}
\usepackage{amssymb,amsmath}
\usepackage{graphicx, subfigure, wrapfig}
\usepackage{amsfonts}
\usepackage[english]{babel}

\def\0{{\bf 0}}

\def\r{{\bf r}}
\def\p{{\mathbf{p}}}
\def\q{{\mathbf{q}}}
\def\r{{\mathbf{r}}}
\def\s{{\mathbf{s}}}

\def\z{{\mathbf{z}}}

\def\Sigma{\boldsymbol{\sigma}}

\newtheorem{theorem}{Theorem}[section]
\newtheorem{proposition}[theorem]{Proposition}

\newtheorem{remark}[theorem]{Remark}

\newtheorem{definition}[theorem]{Definition}

\begin{document}

\title{Notes on relative equilibria of isosceles  molecules in classical approximation}

\author{Damaris McKinley\footnote{Wilfrid Laurier University, 75 Univ. Av. West, Waterloo, N2L 3C5, Ontario, Canada. Email: damaris.duma@gmail.com} \,\,
Daniel Pa\c sca\footnote{Department of Mathematics and Informatics, University of Oradea, University Street 1, 410087, Oradea, Romania. Email: dpasca@uoradea.ro}\,\,
Cristina Stoica\footnote{Wilfrid Laurier University, 75 Univ. Av. West, Waterloo, N2L 3C5, Ontario, Canada. Email: cstoica@wlu.ca}}

\maketitle


\noindent
\textbf{Abstract} We study a classical  model  of isosceles triatomic ``A-B-A'' molecules. The atoms, considered mass points,   interact mutually via a  generic  repulsive-attractive binary potential. 
First we show that the  steady states, or relative equilibria (RE), corresponding to rotations about the molecule symmetry axis may be  determined  qualitatively  assuming  the knowledge of 1) the shape of the binary interaction potential, 2) the equilibrium  diatomic distances (i.e., the equilibrium bond length)  of the  A-A and A-B molecules, and 3)  the  distance at which the RE of the diatomic A-A molecule ceases to exist. No analytic expression for the interaction potentials is needed. 
Second we determine  the stability of the isosceles RE  modulo rotations using geometric mechanics methods and  using  Lennard-Jones  diatomic potentials. As a by-product, we verify  the qualitative results on RE existence and bifurcation. 
For isosceles RE we employ the \textit{Reduced Energy-Momentum} method  presented in [J.E. Marsden, \textit{Lectures in Mechanics}, Cambridge University Press, 1992], whereas for   linear (trivial isosceles) RE we introduce   the  \textit{Symplectic Slice} method, a technique based on the findings in the paper [R.M. Roberts, T. Schmah and C. Stoica, \textit{Relative equilibria for systems with configurations space isotropy}, J. Geom. Phys., 56, 762, (2006)].


\section{Introduction}
   In molecular dynamics, classical approximations are used in order to reduce the computational effort required when using the corresponding quantum-based models. These approximations consist in modeling the atoms in the realm of classical mechanics as mass points that interact mutually via a repulsive-attractive potential. The latter is either determined experimentally and tabulated consequently, or it is described via a functional relation with parameters    fitted experimentally.

   A functional relation expressing the atomic repulsive-attractive interaction is the $12$-$6$ Lennard-Jones potential \cite{Lennard-Jones (1881)}.  
Its parameters  are derived from ``second-virial" or ``viscosity coefficients" \cite{BCH}, providing a particularly accurate model for neutral atoms and molecules, as well as noble gas atoms. However, the use of the powers $12$ and $6$  in the potential is due to a computational convenience \cite{Brush70}. For instance,   the power of the repulsive inverse of the distance term is slightly different from one gas to another, taking values in between 12 and 13. 
Thus   a qualitative analysis of an atomic model within the framework of classical mechanics need not be restricted to the $12$-$6$ potential, but rather to a repulsive-attractive  binary interaction  given by some generic graphical representation.

 As discussed in the physical-chemistry literature \cite{Pav, Ko99}, \cite{Ko00},  rotating steady states about stationary axes during which the ``shape" of a molecule does not change  can be used to explain and predict features of quantum spectra. In classical mechanics, such states are  called \textit{relative equilibria} (RE). A detailed study of triatomic  molecules' RE can be found in \cite{Ko99, Ko00}; in these papers   the authors use  experimentally determined potentials. Also, criteria of bifurcations of RE  from equilibria    are determined in \cite{MR99}.

\begin{figure}
\centerline
{\includegraphics[scale=0.9]{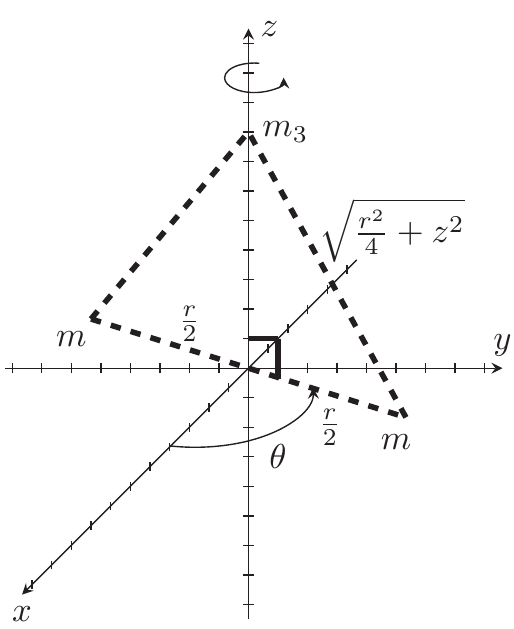}}
{\caption{ \label{threemasspointsystemchgcrds} The model of an isosceles molecules considered here is given by a   three-degrees of freedom system with  coordinates are $r$, $z$ and the angle $\theta$. The molecule is rotated about the $Oz$ axis and (by symmetry) it maintains its isosceles shape at all time. 
}}
\end{figure}

  Here we consider  RE  of ``A-B-A"-type molecules in classical approximation. Such molecules are formed by three atoms that form an isosceles triangle, with two identical A atoms located at the vertices of the  triangle's base,  and the third $B$ atom  at the third vertex (see Figure \ref{threemasspointsystemchgcrds}).  We show  that    the isosceles RE and their bifurcations may be determined qualitatively when no 
analytical expression of the repulsive-attractive molecular diatomic potentials is available. 
The necessary and sufficient information required consists of  three  parameters: the equilibrium  diatomic distances   of the   A-A and A-B molecules  and diatomic distance at which the RE of the A-A molecule ceases to exist (see Proposition \ref{summary_prop}). 
Further, we determine the  isosceles RE stability in the full phase-space using geometric methods, in  particular the Reduced Energy  Momentum (REM) method of Marsden and co-workers \cite{Ma92} and  the \textit{Symplectic Slice} method based on the findings in \cite{RSS06}. We implement these on some numerical examples using Lennard-Jones-type potentials, which also verify the qualitative results on the RE existence and bifurcations.

  While the existence of  isosceles RE  is  qualitatively determined    using three parameters only,   establishing  RE stability  require   additional ``second order" parameters that lead to cumbersome qualitative classifications diagrams, beyond our scope. For example, one needs to know the \textit{dissociation} distance, that is the distance at which   the diatomic potential changes its concavity.   
   Thus we  chose to calculate the RE stability assuming  the knowledge  the inter-atomic A-A and A-B potentials. 
   
   Recall that the stability of an equilibrium calculated on an invariant manifold of a ODE does not imply  stability in the full phase space; in our case, stability within the isosceles invariant manifold may be lost in  directions that breaks the symmetry of the triangle. Thus we must consider the molecule  as a 6 degree of freedom mechanical system (using Jacobi coordinates) with spatial rotational symmetry and apply appropriate stability methods that ignore the neutral symmetry directions.
N-body systems with spatial rotational symmetry stability involve    the moment of inertia tensor  at the RE configuration, represented by a $3 \times 3$ matrix and  its inverse. For isosceles RE the inertia matrix is non-singular, and in this case  we employ the phase-space splitting of the Reduced Energy-Momentum  (REM) method of Marsden and co-authors \cite{Ma92}. 
 For linear RE, the inertia matrix is singular (having a zero eigenvalue), and thus the splittings of the  REM method are not applicable. In this case we introduce   the \textit{Symplectic Slice} method,  based on the theoretical work in \cite{RSS06} (but also see \cite{Sc11}). Practically, we apply a  change of coordinates  that re-organises the dynamics as a coupled system of a rigid-body-like system (corresponding to the invertible part of the locked inertia tensor $\mathbb{I}$) and a  simple mechanical system of type "kinetic+ potential". Note that while the theory was formulated previously, to our knowledge,  it is the first time the Symplectic Slice method is implemented on a concrete example. Future work will extend its use to more complex systems.

We  verify our qualitative findings in RE existence and bifurcation and calculate stability  
numerically using  Lennard-Jones-type binary potentials \cite{Lennard-Jones (1881)}.
 For interactions between atoms of different species we adopt the Lorentz-Berthelot rule (\cite{{Lorentz81}, {Berthelot89}, {Kirchner12}}), a rule that assumes the  knowledge of   the coefficients of identical atoms. 
We consider   examples  in the  cases  when the outer atoms have greater, equal and smaller mass then the central atom; we call these models $D_2H$, $H_3$ and $H_2D,$ respectively. For each case we present the RE \textit{Energy-Momentum} (EM) diagram. 

 The paper is organized as follows: first we define    binary interactions  qualitatively  and  provide a generic picture of the RE of  diatomic molecules in classical approximation.   In Section  \ref{Ch_Mod} we introduce  isosceles triatomic molecules and determine the RE  equations. In Section \ref{sect:_RE} we establish criteria for the existence and bifurcations of the RE. In Section \ref{sect:stab} we describe the REM and the Symplectic Slice methods. In the next Section we consider specific Lennard-Jones-type binary  potentials and calculate numerically  the EM diagrams. We end with Conclusions.

\section{A generic model for molecular interactions}
 \label{sect:pot_generic}

\begin{figure}
\center
\includegraphics[scale=0.4]{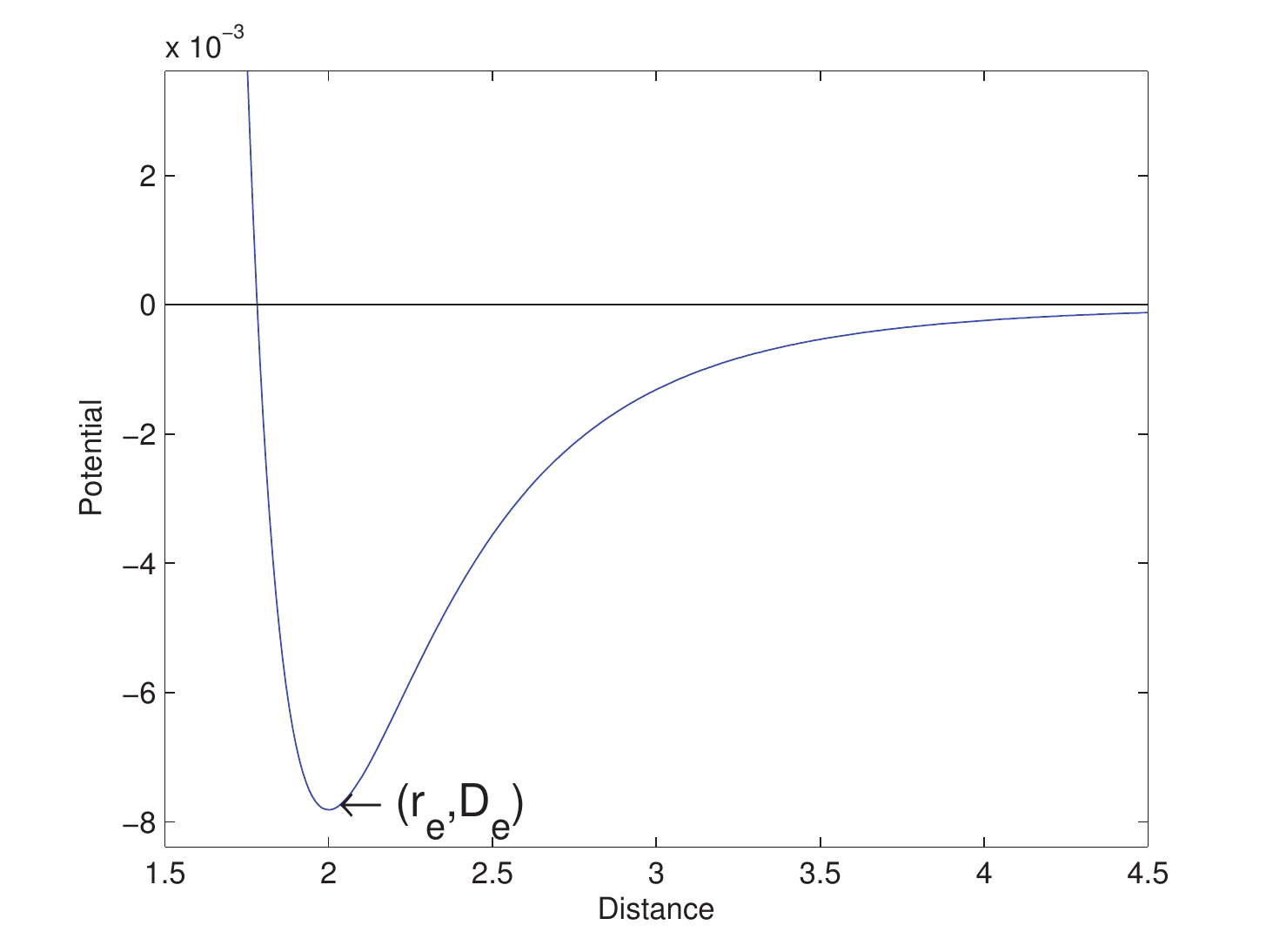}
{\caption{A typical shape of a repulsive-attractive inter-atomic potential. The parameters $r_e$ and $D_e$ are the equilibrium diatomic distance, and the depth of the potential well, respectively.}}
\label{generalpotential}
\end{figure}

Classical approximations of molecules  consist in $N$-body problems with mutual binary interaction 
given by a repulsive-attractive smooth potential as  in Figure  \ref{generalpotential}.
 The atoms are idealized  mass points,  and the interatomic potential  is either determined directly from experiments, or is described via a functional relation (for example, a $12$-$6$ Lennard-Jones potential) with parameters  fitted experimentally. 

\medskip
\noindent
\begin{definition}  \label{assump_0}  { \bf \text{A  generic binary molecular potential}} is a function $U: (0, \infty) \to \mathbb{R}$ such that 

\begin{enumerate}

\item  $\displaystyle{ \lim\limits_{r \to 0^+} U(r)=\infty}\,$ (the atoms repel when close);

\item  the function $U$ has a unique minimum, and no other critical point (the atoms  attract each other at medium range, displaying a  unique  stable equilibrium position);
 
\item $\displaystyle{ \lim\limits_{r \to \infty} U(r)=0}\,$ (at long range     the interaction vanishes);

\item $U(r)$ is  $o\left(r^2\right)\,$ (the binary long-range interaction decays faster that  the centrifugal forces). 

\end{enumerate}

\end{definition}

\begin{definition} \label{assump_1} { \bf \text{A  molecule in classical approximation}} consists of $N$ mass points interacting mutually via a  generic binary molecular potential. 

\end{definition}

\begin{figure}[h!]
\centering
   \includegraphics[scale=0.31]{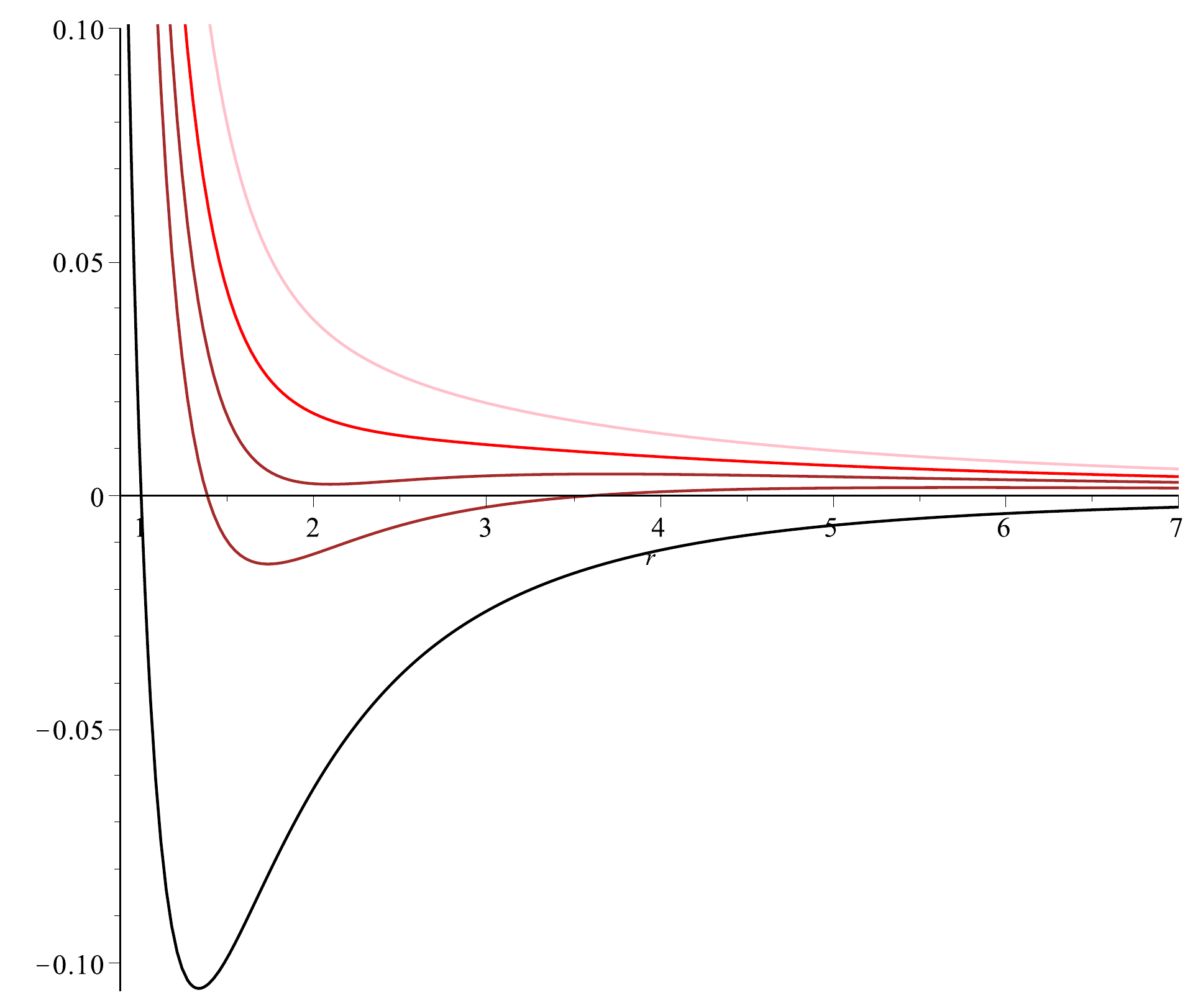}
   \caption{The amended potential $U_{c}(r)$ for various values of $c$. The black curve corresponds to $c=0$ (color online).}
      \label{di_RE_plot}
  \end{figure}

\medskip
\noindent
Consider  the classical approximation of a molecule  formed by two atoms $A$  and $B$ of mass $m_A$ and $m_B$, respectively. The molecule  is  rotated with constant angular velocity about an axis perpendicular to the centre of mass of the $AB$ system.  The  rotating steady states, or relative equilibria (RE), are the states where  the (inter-atomic)  attractive and centrifugal forces  are balanced. The distance between the atoms at a RE is determined as   a critical point of the \textit{amended (or reduced, or effective) potential}
\[U_{c}(r):=\frac{c^2}{mr^2}+U(r)\,,\]
where $c$ denotes the angular momentum and $m$ is the relative mass $m=m_Am_B /(m_A+m_B)$.

In Figure \ref{di_RE_plot} we sketch a  plot of the augmented  potential for various momentum values. The  relative equilibria  (RE)  are found  as the critical points of $U_{c}(r)$. At momentum zero (black curve), apart from the equilibrium configuration  (at bottom of the well), a fictitious equilibrium is located at $r=\infty$. 
For small non-zero momenta (brown curves) there are two RE: one stable (minimum),  emerging from the equilibrium at the bottom of the well, and  one unstable (maximum), emerging from the fictitious equilibrium at infinity. As the momentum increases, at critical value $c_{\text{0}},$ the two RE merge and  cease to exist  at a cusp. The red and pink curves correspond to $c_{\text{0}}$ and  a value $c> c_{0}$, respectively. 

\begin{figure}[h!]
\centering
              \includegraphics[scale=0.34]{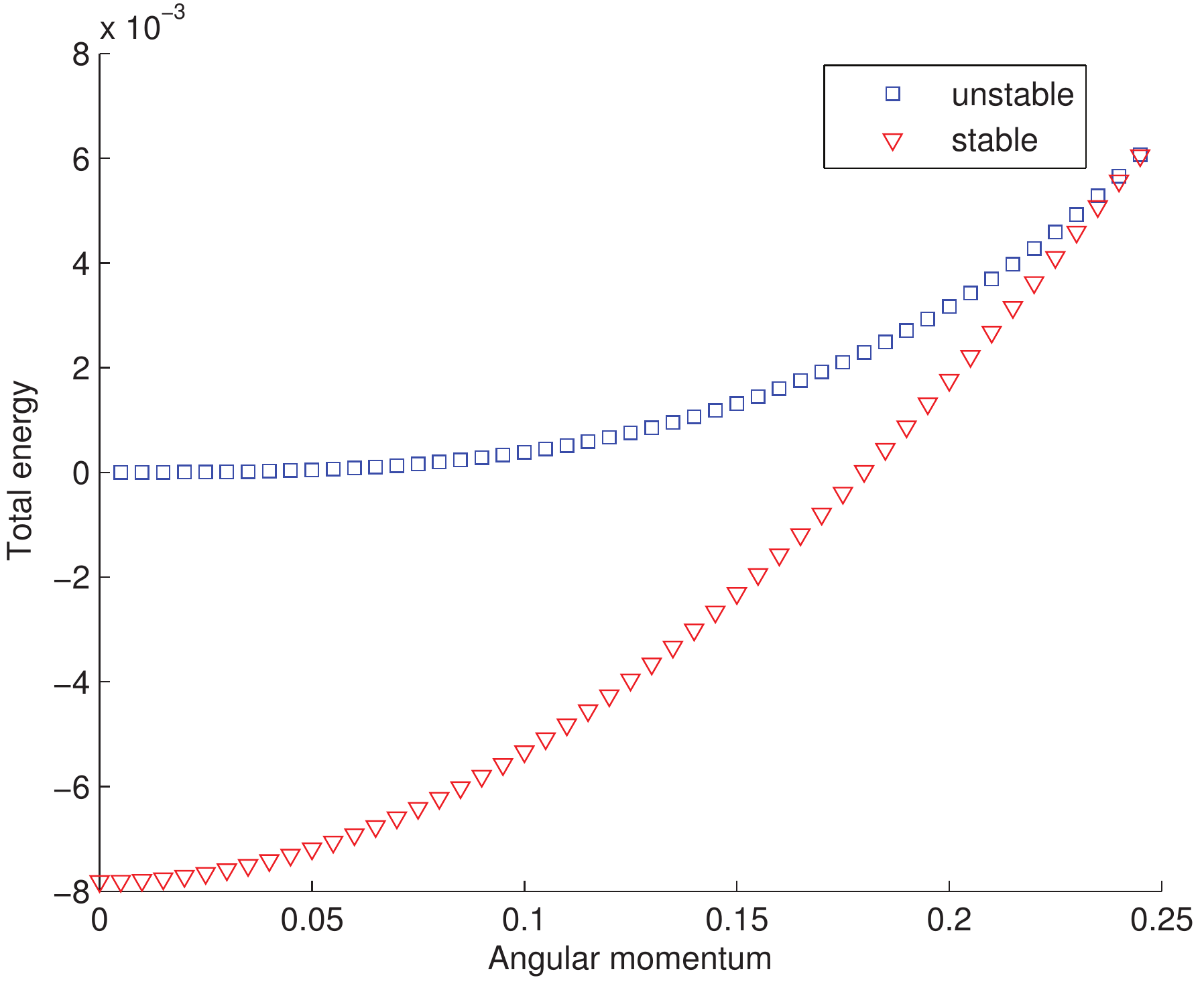}
     \caption{A EM map for a generic diatomic molecule (color online).}
     \label{EM_dia}
  \end{figure}

Figure \ref{EM_dia} depicts  a  energy-momentum (EM) diagram.
Each point represents a RE in  coordinates given by angular momentum against total  energy,   At zero momentum we have the equilibrium configuration (the minimum of the potential) and a fictitious RE at infinity. 

\section{Isosceles three point-mass systems}\label{Ch_Mod}

\indent
Consider a triatomic molecule  of type A-B-A as a classical three point-mass problem  with  binary repulsive-attractive  interaction.  We set up the dynamics in Jacobi coordinates coordinates $(\r, \s),$  with $\r$ is the relative vector between  the two A atoms,  and $\s$ is the vector between the centre of mass (and also midpoint) of the A atoms system and  the B atom.  The configuration of the atoms is $Q:=\left(\mathbb{R}^{6} \setminus \{collisions\} \right)$  and the dynamics is given by  the Lagrangian system with  $L: TQ \to \mathbb{R}$,

\begin{equation}
L(\boldsymbol{r}, \boldsymbol{s}, \dot {\boldsymbol{r}}, \dot {\boldsymbol{s}})= \frac{M_1}{2} \dot{\boldsymbol{r}}^2 +\frac{M_2}{2}\dot{\boldsymbol{s}}^2  -V(\r, \s)
\label{Lag_3b}
\end{equation}
with $\r=(r_x, r_y, r_z),$ $\s=(s_x, s_y, s_z),$
\begin{equation}
V(\r, \s) = V_{AA}(|\boldsymbol{r}|)+V_{AB} \left(\left|\boldsymbol{s}+(1/2)\boldsymbol{r} \right| \right)+V_{AB} \left(\left|\boldsymbol{s}-(1/2)\boldsymbol{r}\right| \right)
\end{equation}
$M_1=m_A/2$ and $M_2=(2m_Am_B)/(2 m_A+m_B), $ and 
$V_{AA}$\, and $V_{AB}$  the potentials between two A atoms, and A and B atoms, respectively.%

\bigskip
\noindent
Given that two of the mass points are equal,  one can prove  that motions with the three atoms in an isosceles shape at all times form a 6 dimensional invariant submanifold of the flow in the full 12 dimensional phase-space.
 Let the configuration space of the isosceles sub-system  be described  by  the coordinates $(r_x, r_y,s_z)$. Passing to cylindrical coordinates $(r_x, r_y,s_z) \to (r, \theta, z )$  and  using Legendre transform, we obtain the Hamiltonian formulation of the isosceles problem, with a Hamiltonian given by 
\begin{equation}
\label{Hamiltonian_Isosceles}
H(r,\theta, z, p_r, p_\theta, p_z):
=\frac{1}{m_A}\left(p_r^2+\frac{p_\theta^2}{r^2}\right)+\frac{2m_A+m_B}{2m_Am_B}\,p_z^2+F(r)+2G\left(\sqrt{\frac{r^2}{4}+z^2}\right)
\end{equation}
where we  denote  $F:=V_{AA}\,$, $G:=V_{AB}\,.$  From the general theory or by direct verification, we have that the conservation of the energy of the system, and  also,  the conservation of angular momentum $\dot{p}(t)=const.=:c$. Substituting the latter into \eqref{Hamiltonian_Isosceles} we obtain the \textit{reduced Hamiltonian}\\
\begin{equation}
\label{H_red_v11}
H_{\text{red}}(r, z, p_r, p_z; c)=\frac{1}{m_A} \, p_r^2+ \frac{2m_A + m_B}{4m_A m_B}\,p_z^2 + \frac{c^2}{m_Ar^2}+F(r)+2G\left(\sqrt{\frac{r^2}{4}+z^2}\right) \,.
\end{equation}
Thus the dynamics is  reduced to a two-degrees of freedom Hamiltonian system  parametrized by the angular momenta $c$. The equations of motion are
\begin{align}
&\dot{r}=\frac{2}{m_A}\,p_r\,,  \, \quad \quad \quad \quad \quad \dot{p}_r=\frac{2c^2}{m_Ar^3}-F'(r)-G' \left(\sqrt{\frac{r^2}{4}+z^2}\right)\frac{r}{2\sqrt{\frac{r^2}{4}+z^2}}\,, \label{H_red_mod_eq_3}\\
&\dot{z}=\frac{2m_A + m_B}{2m_A m_B}\,p_z\,,
\quad \quad \dot{p}_z=-2G' \left(\sqrt{\frac{r^2}{4}+z^2}\right)\,\frac{z}{\sqrt{\frac{r^2}{4}+z^2}}\,. \label{H_red_mod_eq_4}
\end{align}%

The RE are determined as the equilibria of $H_{\text{red}}$. Note that the equilibria of the un-reduced Hamiltonian \eqref{Hamiltonian_Isosceles}  correspond to RE with zero angular momentum. 
  The RE have $p_r=p_z=0$ and their coordinates  are found by solving the system
\begin{align}
& \frac{2c^2}{m_Ar^3}-F'(r)-G' \left(\sqrt{\frac{r^2}{4}+z^2}\right)\frac{r}{2\sqrt{\frac{r^2}{4}+z^2}}=0\,, \label{systequil_mod_1}\\
&2G' \left(\sqrt{\frac{r^2}{4}+z^2}\right)\,\frac{z}{\sqrt{\frac{r^2}{4}+z^2}} =0\,. \label{systequil_mod_2}
\end{align}
To simplify the presentation, we call RE the solutions of the algebraic system \eqref{systequil_mod_1} - \eqref{systequil_mod_2}.

\bigskip

\begin{remark}
\label{third_mass_remark}
Since the mass parameter $m_B$ does not appear in the system \eqref{systequil_mod_1} - \eqref{systequil_mod_2} the RE do not depend  on the mass $m_B$. (Physically, this is because $m_B$ does not add to the rotational inertia of the system.)
\end{remark}

\section{Relative equilibria}
\label{sect:_RE}

We distinguish 2 classes of RE (modulo the sign of $c$):
1) \textbf{linear} RE, that is  solutions with $z=0$.  In this case the atoms form a rigid steadily rotating segment; and 2)  \textbf{isosceles} RE, that is solutions with  $z\neq 0$. In this case the atoms form a rigid steadily rotating isosceles triangle. Without loosing generality,  we consider $c>0$.

\smallskip
 Recall that the potential $F$ describes the attractive-repulsive interaction between  the identical A atoms, whereas  $G$ describes the attractive-repulsive interaction  between the A and B atoms.
Let the equilibrium diatomic distances be $r_{e}^F$ and $r_{e}^G$
 for the potentials $F$ and $G$, respectively. .

\subsection{Linear RE ($z=0$)}
\label{Planar_equilateral_section}
Let $c$ be fixed. From system (\ref{systequil_mod_1}) - (\ref{systequil_mod_2}) and using that $z=0$,  we have that linear RE are solutions of
$\displaystyle{
(2c^2)/(m_A\,r^3)=F'(r)+G' \left(r/2 \right)\,.
}$
We determine these by   counting of the number of intersections of the curves $\displaystyle{2c^2/(m_Ar^3)}$ and $\displaystyle{F'(r)+G' \left(r/2 \right)}$; each intersection of these curves corresponds to a root of the given equation.

\medskip
\noindent
\textbf{\underline{Equilibria}} ($c=0$) We solve
\begin{equation}
\label{planar_equilibrium}
F'(r)+ G' \left(\frac{r}{2} \right) =0\,
\end{equation}
Given the shape of $F(r)$ and $G(r)$ (see Figure \ref{generalpotential}),  this equation 
can have  one, two or three roots. Assuming that the  function $F(r)+ G \left(\frac{r}{2} \right)$  behaves as a  potential described by  Definition \ref{assump_0}, we consider  that any three atoms in a classical interaction display a unique linear equilibrium configuration and so Equation (\ref{planar_equilibrium}) has a unique solution.

\medskip
\noindent
\textbf{\underline{Relative Equilibria}} ($c\neq0$)  Since the sum $F(r)+ G \left(\frac{r}{2} \right)$ is a generic molecular potential,  linear RE are determined in the same manner as for a  diatomic molecule  (see Section \ref{sect:pot_generic}).

\subsection{Isosceles RE ($z\neq0$)}
For $z \neq 0$ the system ({\ref{systequil_mod_1}}) - ({\ref{systequil_mod_2}}) becomes
\begin{align}
\label{spatial_RE_1}
&\frac{2c^2}{m_A r^3}-F'(r)=0\\
\label{spatial_RE_2}
&G' \left(\sqrt{\frac{r^2}{4}+z^2}\right)=0
\end{align}

\bigskip
\noindent
\textbf{\underline{Equilibria}} $(c=0)$ In this case Equation (\ref{spatial_RE_1}) reduces to $F'(r)=0$, which is solved by $r=r_{e}^F$. Then using  Equation (\ref{spatial_RE_2}),
and the fact that $G'(r_{e}^G) =0$, we obtain
\begin{equation*}\sqrt{\frac{\left(r_{e}^F \right)^2}{4}+z^2} = r_{e}^G \,.
\end{equation*}
Provided $r_{e}^F\leq 2r_{e}^G$
it follows that there are two isosceles equilibria located at
$\displaystyle{\left(r_{e}^F,\,\pm\, \sqrt{\left(r_{e}^G \right)^2  -\frac{ \left(r_{e}^F \right)^2}{4}}\, \right)}$
that  degenerate to a linear RE when $r_{e}^F=2r_{e}^G.$

\begin{figure}
\centerline
{\includegraphics[scale=0.6]{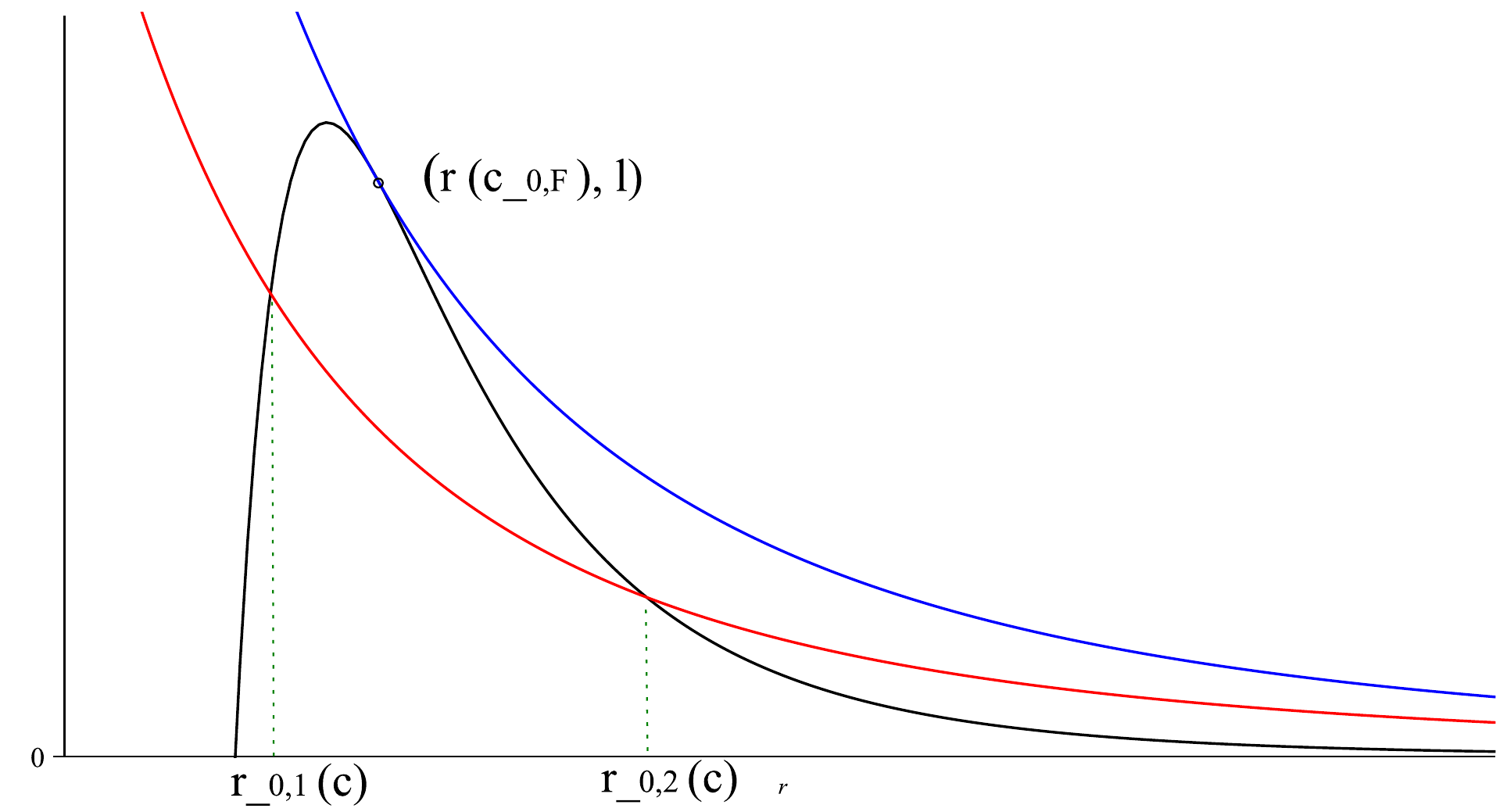}}
{\caption{ \label{qual_spat_RE}  The RE of the diatomic molecule  formed by the identical $A$ atoms   are determined as intersections of $y(r)=F'(r)$ (in black) with curves $\displaystyle{y(r)=c^2/(mr^3)}$ for various values of the angular momentum $c$ (color online). Let $c_{0,F}$ be the critical momentum value above which such RE disappear. For a fixed $c\in (0, c_{0,F})$, we plot   the curve $\displaystyle{y(r)=c^2/(mr^3)}$  (in red) and  denote by $r_{0,1}$ and $r_{0,2}$  the corresponding RE.  
At $c=c_{0,F}$ we have  the curve $\displaystyle{y(r)=c_{0,F}^2/(mr^3)}$ (in blue) and  denote  $\displaystyle{l= r(c_{0,F})}$ the corresponding RE distance.}}
\end{figure}

\bigskip
\noindent
\textbf{\underline{Relative Equilibria}}  $(c\neq 0)$ We observe that the real solutions of the Equation (\ref{spatial_RE_1}) are in fact the critical points of  the amended potential
\begin{equation}
F_{c}(r):= \frac{c^2}{m_Ar^2} + F(r)
\label{F_red}
\end{equation}
associated to the diatomic molecule formed by the A atoms. Thus  isosceles RE exist only for momenta $c$  that allow the existence of RE for the diatomic A-A molecule. 
We denote by $c_{0,F}$  the critical momentum value above which the diatomic  A-A RE disappear, and   by $l$ the corresponding inter-atomic distance, i.e. $l:=r(c_{0,F})$  (see  Figure \ref{qual_spat_RE}).
\noindent
For $0<c\leq c_{0,F},$  let the solutions of Equation (\ref{spatial_RE_1}) be $r_{0,1}(c)$ and  $r_{0,2}(c)$. We have %
$
r_{e}^F< r_{0,1}(c)\leq l \leq r_{0,2}(c)\,.
$
Note that $\lim \limits_{c\to 0} r_{0,1}(c) =r_{e}^F\,,$  $\lim\limits_{c\to 0} r_{0,2}(c)
=\infty,$ and $\lim\limits_{c\to {c_{0,F}}} r_{0,1}(c)  = \lim\limits_{c\to {c_{0,F}}} r_{0,2}(c) = l\,$.
Equation (\ref{spatial_RE_2}) becomes
\begin{equation}
\label{eq:RE_z_spatial}
z^2= \left(r_{e}^G\right)^2 - \frac{r_{0,i}^2(c)} {4}
\end{equation}
and it accepts real  solutions only if
\begin{align}
\label{cond_RE_z_spatial}
r_{0,1}(c) \leq 2r_{e}^G \,\,\,\,\text{or} \,\,\,\,  r_{0,2}(c)\,\leq 2r_{e}^G.
\end{align}
We are able now to distinguish:

\smallskip
\noindent
1.   If $\displaystyle{2r_{e}^G \leq r_{e}^F }\,$, then there are no solutions. %

\smallskip
\noindent
2. 
If $\displaystyle{r_{e}^F < 2r_{e}^G}$ then consider   the momentum value $c_{1,s}$ for which $2r_{e}^G$ is a RE for a diatomic molecule formed  by the identical atoms, that is, $c_{1,s}$ is the momentum value for which 
\begin{equation}
\label{c_2}
F'(2r_{e}^G)=\frac{2c_{1,s}^2}{m_A\left(2r_{e}^G\right)^3}\,.
\end{equation}
We have the following subcases:

\smallskip
\noindent
(a) If $2r_{e}^G <l$ then for $c \in (0, c_{1,s})$  we have one RE family with
\begin{equation}
\label{z_1}
z_{1,2}= \pm \sqrt{\left(r_{e}^G \right)^2  -\frac{ \left(r_{0,1}(c) \right)^2}{4}}\,.
\end{equation}
This family  ceases to exist at  $c=c_{1,s}$ where  it meets a linear family. 

\smallskip
\noindent
(b)  If $\displaystyle{l< 2r_{e}^G}$ then

\smallskip
 (i) for $c \in (0, c_{1,s})$ we have  
one RE family  with
\begin{equation}
\label{z_2}
z_{1,2} = \pm \sqrt{\left(r_{e}^G \right)^2  -\frac{ \left(r_{0,1} \right)^2}{4}}\,;
\end{equation}

\smallskip
 (ii) at $c=c_{1,s},$ a second family ``$z_{3,4}$" of spatial RE branches from a linear RE;

\smallskip
 (iii)  for $c \in(c_{1,s}, c_{0,F})$ we have two RE families with
\begin{equation}
\label{z_2_v2}
z_{1,2} = \pm \sqrt{\left(r_{e}^G \right)^2  -\frac{ \left(r_{0,1} \right)^2}{4}}\,, \quad \quad z_{3,4} = \pm \sqrt{\left(r_{e}^G \right)^2  -\frac{ \left(r_{0,2} \right)^2}{4}}\,;
\end{equation}

\smallskip
 (iv) at $c=c_{0,F}$, the isosceles families ``$z_{1,2}$" and ``$z_{3,4}$" join. No isosceles RE exist for $c>c_{0,F}.$

\noindent
We summarize our findings in:

\begin{proposition}
\label{summary_prop}
Let $r_e^F$ be the equilibrium distance of  the diatomic molecule formed by  the identical A atoms and  with  potential $F$ and let $r_e^G$ be the equilibrium distance of  the diatomic molecule formed by the non-identical atoms A and B   with  potential $G$.
Further, let $c_{0,F}$ be the threshold value of the momentum where the two RE families of the diatomic A-A molecule 
merge and disappear, and  
denote by $l$ the corresponding distance between the atoms; in other words,  let $l:=r(c_{0,F})$. 
Then we have the following:

 \begin{enumerate}

 \item\label{no_spat} If $2 r_{e}^G \leq r_{e}^F$ then there are no isosceles RE.

\item If $r_{e}^F < 2 r_{e}^G$, let us denote by 
$c_{1,s}$  the momentum value for which $2r_{e}^G$ is a RE for a diatomic molecule formed  by the identical atoms and  with  potential $F$.
Then we distinguish:

\begin{enumerate}
\item\label{one_spat} 
If $ 2 r_{e}^G \leq l$, then there is one isosceles family of RE parametrized by $c \in (0, c_{1,s})$. This family emerges from the isosceles equilibrium and  at $c=c_{1,s}$ it joins a linear  family of RE.

\item\label{two_spat}  If $l < 2 r_{e}^G$ then:

\begin{enumerate}
\item  For values $c \in (0, c_{1,s})$ there is one isosceles family of RE which emerges from the isosceles equilibrium.

\item  At $c= c_{1,s}$  there is a bifurcation.  A second isosceles RE family occurs and  the two isosceles families of RE co-exist  for $c \in (c_{1,s}, c_{0,F})$; they merge at $c=c_{0,F}$.

\item For $c>c_{0,F}$, no isosceles RE exist.
 \end{enumerate}

\end{enumerate}

 \end{enumerate}

\end{proposition}

\begin{remark} [Three identical atoms]
\label{remark_3_ident_atoms}
If the three atoms are identical, then $F(r)=G(r)$, and so $r_{e}^F=r_{e}^G.$ In particular, $r_{e}^F < 2r_{e}^G= 2r_{e}^F$ and so  Case 2. of the Proposition above applies.
\end{remark}

\section{Stability}
\label{sect:stab}

We now present two  geometric methods appropriate for computing stability of the isosceles and linear RE, respectively, modulo spatial rotations.  
These methods will be employed on  the numerical examples in the next section.

\smallskip
Broadly speaking, the symmetry reduced Hamiltonian of a mechanical system with rotational symmetry takes the form  ``kinetic plus reduced  potential".  A RE is  stable modulo rotations if the Hessian matrix $D_2V_{\text{amended}} (\r_e,\z_e;c)$ of the amended potential     is positive definite. 
  This test is in fact a specialization of the non-linear (Lyapunov) stability test to mechanical systems with symmetry. The positive definiteness of the amended potential $D_2V_{c} (\r_e,\z_e;c)$ insures that the Hessian of the reduced Hamiltonian $D_2H_{\text{red}}$  is also positive definite and so $(\r_e,\z_e;c)$ corresponds to  a local minimum of the energy. 
 If $D_2V_{c} (\r_e,\z_e;c)$ is not positive definite, then one proceeds to the spectral analysis of the linearization of the vector field. In the canonical Hamiltonian context, linear stability is guaranteed if the linearization matrix is semi-simple and all its eigenvalues have zero real part. 
   Linear  stability does not  ensure stability, but only the existence of  the  local  normal modes (simple harmonic motions) near RE.  However,  linear instability predicts  instability.

\medskip

 In the present  stability calculations   we employ two  methods: for the isosceles RE, we use the splitting of the phase space  of the REM method   (see \cite{Ma92}), whereas for the linear RE we introduce the Symplectic Slice method, that is a method with  theoretical justification in \cite{RSS06} (also, see \cite{Sc11}).     We sketch below the application of these methods to our class  of RE. 
  
\subsection{The REM method for  isosceles RE}

Using the notation from Section \ref{Ch_Mod}, the dynamics of a three mass point system is given by the $SO(3)$-invariant  Hamiltonian 
\[H(\r,\s, \p_\r, \p_\s) = \frac{\p_{\r}^2}{2M_1} +\frac{\p_{\s}^2}{2M_2}+ V(\r,\s)\,.\]
where $(\p_\r, \p_\s)\in\mathbb{R}^3 \times \mathbb{R}^3$ are the generalized momenta corresponding to the coordinates $(\r, \s).$ We denote the Lie algebra and co-Lie algebra of $SO(3)$ by $so(3)$ and $so^*(3)$. 
The (angular) momentum map
corresponding to the action of $SO(3)$ on the phase space is 
$J:T^*Q\rightarrow  so(3)^*$, 
$
\label{mom_here}
J(\r,\s, \p_\r, \p_\s)  = \r \times \p_\r +\s \times \p_\s\,.
$

\smallskip
Consider a  RE  
with angular velocity $\boldsymbol{\xi}_e \in  so(3)$,  angular momentum $\boldsymbol{\mu}_e \in so(3)^*,$ and  base-point  $\q_e=(\r_e, \s_e).$ The REM method uses the augmented Hamiltonian 
 \[H_{{\boldsymbol{\xi}}}(\r,\s, \p_\r, \p_\s ):= H(\r,\s, \p_\r, \p_\s ) - \left<  J(\r,\s, \p_\r, \p_\s )\,, \boldsymbol{\xi} \right>\,,\quad \boldsymbol{\xi} \in so(3)\] 
as a Lyapunov function on  an appropriate splitting of the co-tangent bundle phase-space  variations that has the additional property of bringing the linearized equations of motion into a normal form. The stability computation  is reduced to a test of positive definiteness of the second variation $\delta^2V_{\boldsymbol{\mu}_e}$ at  $\q_e$   of the amended potential 
\[V_{\boldsymbol{\mu}_e} (\r, \s) :=  V(\r, \s) + \frac{1}{2} \left< \boldsymbol{\mu}_e, \mathbb{I}^{-1}(\r, \s) \boldsymbol{\mu}_e \right>\]
restricted to a subspace of configuration variations. 
Briefly,  the tangent to the  configuration space at the RE  allows a three-way splitting $T_{\q_e}Q= {\cal M} \oplus \mathcal{V}_{RIG} \oplus \mathcal{V}_{INT}$. 
The   subspace ${\cal M}$  contains variations  along the symmetry group and  the   corresponding $\displaystyle{\left[\delta^2V_{\mu_e}\big|_{\cal M} \right]}$ matrix block is always zero.
  The subspace $\mathcal{V}_{RIG}$ contains \textit{rigid} variations that, loosely speaking,  do not affect the shape of the  RE configuration which is  perceived  as a rigid body. 
 The subspace $\mathcal{V}_{INT}$ is normal to ${\cal M} \oplus \mathcal{V}_{RIG}$ and is formed by  \textit{internal} variations which deform the RE configuration.  The  REM phase-space splittings   are explained in detail in Section 2  of \cite{ScSt06}.

\smallskip
For  isosceles RE, the variations along ${\cal M}$ are given by infinitesimal rotations of the  triangle  along the symmetry axis;  the variations along $\mathcal{V}_{RIG}$ take the ``rigid" RE triangle out of its plane; and the variations  along $\mathcal{V}_{INT}$   change   the shape the RE triangle.  
Isosceles   RE are of  the form $\q_e=(\r_e, \s_e)$, with $\r_e=(0, r_e,0)$ and $\s_e=(0,0,s_e),$
and with angular velocity $\boldsymbol{\xi}_e = (0,0,\xi_e) \in \mathbb{R}^3 \simeq so(3)$  and angular momentum $\boldsymbol{\mu}_e = (0,0,\mu_e)\in \mathbb{R}^3 \simeq so(3)^*.$ The moment of inertia at a RE is
\begin{equation}
\mathbb{I}(\r_e, \s_e) = 
\left[
\begin{array}{ccc}
M_1r_e^2 + M_2 s_e^2 &0&0\\
0& M_2 s_e^2 & 0\\
0&0& M_1r_e^2\,
\end{array}
\right]\,.
\end{equation}
It is clear that $\mathbb{I}(\r_e, \s_e)$ becomes singular at linear configurations (where $s_e=0$) and that  in this case, one needs a modified method to calculate stability (treated in the next subsection). 

\smallskip
Assuming $\s_e \neq 0$, 
denoting 
${so(3)}_z^\perp:= 
\{ \boldsymbol{\omega} \in so(3)\,|\,  \boldsymbol{\omega} =  (\omega_x, \omega_y, 0)\}$ and  using the REM block diagonalisation 
of $\displaystyle{\delta^2V_{\boldsymbol{\mu}_e} \big|_{\mathcal{V}_{INT}}}$, after  calculations, we have that an isosceles RE is Lyapunov stable if 
\begin{enumerate}

\item the \textit{Arnold} two-form ${\cal A }_{\boldsymbol{\mu}_e} : {so(3)}_z^\perp \times {so(3)}_z^\perp \to \mathbb{R}$, which has the matrix
\begin{align}
\left[{\cal A }_{\boldsymbol{\mu}_e} \right] =  \mu_e^2 \left[\mathbb{I}^{-1}(\r_e, \s_e) - \lambda_e^{-1} \text{I}_3 \right]= \mu_e^2
\left[
\begin{array}{cc}
 \frac{1}{M_2 s_e^2 } - \frac{1}{M_1r_e^2}  &0\\
0 & \frac{1}{M_1r_e^2+ M_2 s_e^2 } - \frac{1}{M_1r_e^2}
\end{array}
\right]
\end{align}
is positive definite, and

\item the second variation $\delta^2V_{\boldsymbol{\mu}_e}$  is positive definite on $\mathcal{V}_{INT}=\{(\delta \r, \delta \s) \in \mathbb{R}^3 \times \mathbb{R}^3 \,|\,   \delta r_x= \delta s_x =0\,, M_1 r_e\,\delta r_z +  M_2 s_e \,\delta s_y =0\}$. Choosing the parametrization $(\delta r_y, \delta s_y, \delta s_z) \in \mathbb{R}^3$ for $\mathcal{V}_{INT}$, we calculate

\begin{align*}
\label{stability-matrix}
&
\delta^2V_{\boldsymbol{\mu}_e} \big|_{\mathcal{V}_{INT}}
=  \left[ {
\begin{array}{ccc} V_{r_y r_y} & \alpha V_{r_y r_z}+V_{r_y s_y}& V_{r_y s_z}
\\
\alpha V_{r_y r_z}+V_{r_y s_y} & 2 \alpha  V_{r_z s_y} +  V_{s_y s_y}  & \alpha  V_{r_z  s_z}+ V_{s_y s_z}
\\ V_{r_y s_z} &
\alpha  V_{r_z  s_z}+ V_{s_y s_z}
& V_{s_z s_z}
\end{array}}\right]
 +
\left(\xi_e \right)_z^2 \left[ {
\begin{array}{ccc} M_1 & 0 & 0
\\
0 & 4 M_1-M_2 & 0
\\ 0 & 0
& 0
\end{array}}\right],
\end{align*}
where
$\displaystyle{\alpha:=- \frac{M_2 s_e}{M_1 r_e}}.$

\end{enumerate}

\smallskip
\noindent
From Condition 1. we deduce that the Arnold form is never positive definite,  and so

\begin{remark} \label{non_stab_crit} An isosceles RE rotating about its symmetry axis is never Lyapunov stable. 
\end{remark}
We then  proceed to study  linear stability and compute the eigenvalues of the linearization $L$ of the Hamiltonian vector-field $X_{H_{{\boldsymbol{\xi}}_e}}$ at $(\r_e, \s_e).$ As known, for a canonical Hamiltonian system, if the matrix of $L$ is semi-simple and all its eigenvalues have zero real part, then the RE is linearly stable. 
We calculate
\begin{align}
L= \left[
\begin{array}{ccc}
A & \mathbb{O} &\mathbb{O} \\
\mathbb{O} &\mathbb{O} & B^{-1} \delta^2 K_{\boldsymbol{\mu}_e} \\
\mathbb{O} & -B^{-1} \delta^2V_{\boldsymbol{\mu}_e} \big|_{\mathcal{V}_{INT}} & \mathbb{O} 
\end{array}
\right] 
\end{align}
where 
\begin{align}
A=
\left[
\begin{array}{cc}
0 & -\frac{1}{M_1r_e^2+ M_2 s_e^2 } + \frac{1}{M_1r_e^2}  \\
 \frac{1}{M_2 s_e^2 } - \frac{1}{M_1r_e^2} & 0
\end{array}
\right]
\end{align}
\begin{align}
B=
\left[
\begin{array}{ccc}
1& 0 & 0  \\
0& 1-  \frac{M_2 s_e^2}{M_1 r_e^2 }  & 0 \\
0&0&1
\end{array}
\right] \quad \quad \text{and}\quad \quad \delta^2 K_{\boldsymbol{\mu}_e}=
\left[
\begin{array}{ccc}
\frac{1}{M_1}& 0 & 0  \\
0& \frac{s_e^2}{r_e^2 } \frac{1}{M_1}+ \frac{1}{M_2} & 0 \\
0&0&\frac{1}{M_2}
\end{array}
\right]\,.
\end{align}

\smallskip
\noindent
The eigenvalues of the block $A$ are solutions of 
\begin{align}
\lambda^2 - \frac{M_2s_e^2 (M_1r_e^2-M_2s_e^2)}{M_2s_e^2(M_1r_e^2+M_2s_e^2)(M_1r_e^2)^2} =0\,.
\end{align}
We observe that for 
$
M_1r_e^2>M_2s_e^2$  the real part of the eigenvalues differs from zero. We immediately deduce:
\begin{remark} A necessary condition for linear stability of an isosceles RE is that $M_1r_e^2<M_2s_e^2$, or equivalently, that the spin axis of the RE  must  be its shortest  principal axis of inertia.
\end{remark}
\begin{remark} From the calculations above, we deduce that linear RE may be  (Lyapunov or linear) stable. 
\end{remark}
Finally, we need  to calculate  the eigenvalues of the remaining block
\begin{align}
\left[
\begin{array}{cc}
\mathbb{O} & B^{-1} \delta^2 K_{\boldsymbol{\mu}_e} \\
 -B^{-1} \delta^2V_{\boldsymbol{\mu}_e} \big|_{\mathcal{V}_{INT}} & \mathbb{O} 
\end{array}
\right]\,. 
\end{align}

\bigskip

\subsection{The Symplectic Slice method for   linear RE}

The stability of linear RE is determined using  the theoretical findings in  \cite{RSS06} and \cite{Sc11}) that we briefly present below. 

\medskip
Let $G$ be a Lie group, with Lie algebra $\mathfrak{g},$ and consider
a smooth left action of $G$ on a finite-dimensional manifold $Q,$
written $(g,\q)\mapsto g\cdot \q.$
For every
$\boldsymbol{\xi} \in\mathfrak{g}$  and $\q \in Q,$ the \emph{infinitesimal action} of $\boldsymbol{\xi}$
on $\q$ is $\boldsymbol{\xi}\cdot \q=\frac{d}{dt}\left.\exp\left(t\boldsymbol{\xi}\right)\cdot z\right|_{t=0}.$
The \emph{isotropy subgroup} of $\q\in Q$ is $G_{\q}:=\left\{ g\in G\mid g\cdot \q=\q\right\} .$
For three-body $SO(3)$-symmetric problems, the action of $SO(3)$ on the set of linear configurations has a non-trivial isotropy group given by a copy of $SO(2).$ In our specific case,   linear RE configurations $(\r_e, \s_e)$, $\r_e= (0, r_e, 0)$\,, $\s_e=(0,0,0)$, have isotropy group  $SO(2)_y$, the subgroup of rotations about the $Oy$ axis.

\smallskip
The geometric framework for studying simple mechanical systems with configurations space (continuous) isotropy is given by degenerate parametrisations, called \textit{slice coordinates}, for neighbourhoods of the isotropic points. The configuration space in the neighbourhood of an isotropic   point is modeled as a twisted product, which is the base space of a principal bundle with fiber the isotropy group of the given point.
Denote by $\mathfrak{g}_{\q_e}$ the isotropy Lie algebra  of a point $\q_e$ and  let  $\mathfrak{g}_{\q_e}^{\perp}$ be a $G_{\q_e}$-invariant complement of $\mathfrak{g}_{\q_e}$.
The slice, $S$,   is a vector space normal to $\mathfrak{g}_{\q_e} \cdot \q_e$;    
 locally we have the identification  $Q \simeq G \times S$ with $\q_e \simeq (e,0) \in G \times S$ (where $e$ denotes the identity in $G$). 
 Further,  
$TQ \simeq G \times \mathfrak{g}_{\q_e}^{\perp} \times TS$, $T^*Q \simeq G \times (\mathfrak{g}_{\q_e}^{\perp})^* \times T^*S$, and  the dynamics may be expressed in $ \mathfrak{g}_{\q_e}^{\perp} \times TS$, or  $\left(\mathfrak{g}_{\q_e}^{\perp}\right)^* \times T^*S$, respectively, coordinates. In our context, we have
 $\mathfrak{g}_{\q_e}=so(2)_y$ and $\mathfrak{g}_{\q_e}^{\perp}=so(2)_{x,z}$ (infinitesimal rotations about $Ox$ and $Oz$ axes).

\smallskip
Practically, we apply a  change of coordinates that  re-organises the dynamics as a coupled system with two parts: a rigid-body-like system corresponding to the invertible part of the locked inertia tensor $\mathbb{I}$, and a  simple mechanical system. The metric  matrix  is re-arranged  into a \textit{reduced locked inertia} $\mathbb{I}_{\text{red}}$ block and a 
\textit{reduced mass}  $\mathbb{M_{\text{red}}}$ block, 
the two being coupled by a \textit{Coriolis}-type term 
$\mathbb{C}$. Next we write the Lagrangian in these coordinates,  pass to the  Hamiltonian side and calculate the linearization of the corresponding Hamiltonian system.

\smallskip
We calculate:
\begin{equation}
S=\{ (\r,\s)=(0, r_e+\sigma_1, 0, \sigma_2, \sigma_3, \sigma_4) \,|\, \sigma_i \in \mathbb{R} \} \simeq \{ \Sigma \,|\,\Sigma= (\sigma_1,  \sigma_2, \sigma_3, \sigma_4) \in \mathbb{R}^4 \}
\end{equation}
\begin{equation}
\mathbb{I}_{\text{red}} (\Sigma)= 
\left[
\begin{array}{cc}
M_1(r_e+\sigma_1)^2 + M_2 (\sigma_3^2+\sigma_4^2)& -M_2 \sigma_2\sigma_4\\
-M_2 \sigma_2\sigma_4 & M_1(r_e+\sigma_1)^2 + M_2 (\sigma_3^2+\sigma_3^2)
\end{array} 
\right]\,,
\end{equation}
\begin{equation}
\mathbb{C}(\Sigma)=
\left[
\begin{array}{cccc}
0&0& M_2 \sigma_4&  M_2 \sigma_3  \\
0 & -M_2 \sigma_3 & M_2 \sigma_2 & 0
\end{array} 
\right]\,,
\end{equation}
\begin{equation}
\mathbb{M}(\Sigma)=
\left[
\begin{array}{cccc}
M_1 &0&0 &0  \\
0 &M_2&0 &0  \\
0 &0&M_2 &0  \\
0 &0&0 &M_2  
\end{array} 
\right]=: \mathbb{M}\,.
\end{equation}
Note that the RE $(\r_e, \s_e)$ corresponds to ${\bf 0}= (0,0,0,0)\in S.$ In slice coordinates the  Lagrangian \eqref{Lag_3b} becomes $L: so(3)_{x,z} \times TS  \to \mathbb{R}$
 \begin{equation}
L(\xi_x, \xi_z, \Sigma, \dot \Sigma) =
\frac{1}{2}\left[\xi_x, \xi_z \right] \mathbb{I}_{\text{red}}( \Sigma) \left[
\begin{array}{cc}
\xi_x  \\
\xi_z 
\end{array} 
\right]+ \left[\xi_x, \xi_z \right] \mathbb{C}(\Sigma) \dot \Sigma + \frac{1}{2} \dot \Sigma^T \mathbb{M}( \Sigma)\dot \Sigma - V(\Sigma).
\end{equation}
Applying the Legendre transform, the reduced Hamiltonian takes the form $H: so(3)_{x,z}^* \times T^*S \to \mathbb{R}$
 \begin{align}
H(\mu_x, \mu_z, \Sigma, \p_{\Sigma}) &= 
\frac{1}{2}\left[\mu_x, \mu_z \right] \mathbb{I}_{\text{red}}^{-1}(\Sigma)
 \left[
\begin{array}{cc}
\mu_x  \\
\mu_z 
\end{array} 
\right]  \nonumber \\
&+ 
\frac{1}{2} 
\left( \p_{\Sigma} - \mathbb{A}(\Sigma) \left[
\begin{array}{cc}
\mu_x  \\
\mu_z 
\end{array} 
\right]  \right)^T 
\mathbb{M}^{-1 }
\left( \p_{\Sigma} - \mathbb{A}(\Sigma) \left[
\begin{array}{cc}
\mu_x  \\
\mu_z 
\end{array} 
\right]  \right)
 + V(\Sigma)
\end{align}
where $\mathbb{A}(\Sigma) = \mathbb{I}_{\text{red}}^{-1}(\Sigma) \mathbb{C}(\Sigma)$ and $\mathbb{A}({\bf 0})=0.$ Finally, to decide stability, we  calculate Hessian $D_2H(0, 0, {\bf 0}, {\bf 0})$ and, if the later is  not positive definite, the spectrum of the linearization of corresponding Hamitonian system  at the equilibrium $(0, 0, {\bf 0}, {\bf 0}).$

\section{Numerical examples}
\label{sect:num}

In this section we determine the RE numerically,  verifying the qualitative findings stated in Proposition   \ref{summary_prop} and apply the stability methods  presented in the previous section. Then we determining the stability of RE using the above geometric.

As apparent from the previous section,   determining the stability of RE  involves the second order derivatives of the diatomic potentials at the RE. Specifically, one must consider the relative positions of the  critical points  and the inflection points (called \textit{dissociation} distances in chemical-physics terminology), leading  to cumbersome classification diagrams.
Rather then focusing on such an analysis, we  use explicit Lennard-Jones type models.

We consider  triatomic molecules with outer atoms with greater, equal and smaller mass then the central atom;  we call these  $D_2H$, $H_3$ and $H_2D$-type molecules, respectively.

As mentioned above, the numerical experiments are performed considering  Lennard-Jones potentials. For interactions between atoms of different species we adopt the average rule of Lorentz-Berthelot  (\cite{{Lorentz81}, {Berthelot89}, {Kirchner12}}), a rule that states that if  atoms of species $i$  and $j$ have binary interacting coefficients $a_{ii}, b_{ii}$ and $a_{jj}, b_{jj},$ respectively, then the coefficients of the interacting potential between the atoms $i$  and $j$  are given by
$a_{ij}:=(a_{ii}+a_{jj})/2$\,and  $b_{ij}:=(b_{ii}+b_{jj})/2.$

\smallskip
In our context, this means that  if the binary potential between the identical A atoms is
$
V_{AA}=F(r)=-{a_{11}}/{r^6} + {b_{11}}/{r^{12}}
$
and the binary potential between  the identical B atoms is
$
V_{BB} = \tilde F(r)=-{a_{22}}/{r^6} + {b_{22}}/{r^{12}}
$
then
$
V_{AB}= G(r)=-(a_{11}+a_{22})/(2r^6) + (b_{11}+b_{22})/(2r^{12})\,.
$
Consequently, we calculate
\begin{equation}
r_e^F=2^{1/6} \left(\frac{b_{11}}{a_{11}}\right)^{1/6}\,,\quad l=5^{1/6} \left(\frac{b_{11}}{a_{11}}\right)^{1/6}\,,\quad r_e^G=2^{1/6} \left(\frac{b_{11}+b_{22}}{a_{11}+a_{22}}\right)^{1/6}\,.
\end{equation}

The criteria in Proposition \ref{summary_prop} become:

\begin{enumerate}

\item Condition $r_e^F \geq 2r_e^G$  (no isosceles RE): 
$\displaystyle{
 \,\,\,2^{\,6} \left(1+ \frac{b_{22}}{b_{11}} \right) -1\leq \frac{a_{22}}{a_{11}}
}$

\item Condition $r_e^F< 2r_e^G$: $\displaystyle{\,\,\,\,\frac{a_{22}}{a_{11}}  < 2^{\,6} \left( 1+ \frac{b_{22}}{b_{11}} \right)-1\,.}$

\begin{enumerate}

\item Condition $2r_e^G<l$  (one family of isosceles RE): $\displaystyle{\,\,\,2^7 \left(1+ \frac{b_{22}}{b_{11}} \right) < 5 \left( 1+ \frac{a_{22}}{a_{11}} \right)\,.}$
Together with the previous condition,  this constraints the ratio ${b_{22}}/{b_{11}}$ to
\begin{equation*}
\frac{b_{22}}{b_{11}} \in \left( \frac{1}{2^6} \left(1+ \frac{a_{22}}{A_{11}}\right) -1\,,  \,\frac{5}{2^7} \left(1+\frac{a_{22}}{A_{11}}\right) - 1 \right)
\end{equation*}

\item Condition $l< 2r_e^G$ (two families of isosceles RE):
$\displaystyle{\,\,\,
5 \left( 1+ \frac{a_{22}}{a_{11}} \right)< 2^7 \left(1+ \frac{b_{22}}{b_{11}} \right) \,.
}$
which is equivalent to
$\displaystyle{\,\,\,
\frac{b_{22}}{b_{11}} > \frac{5}{2^7} \left(1+\frac{a_{22}}{a_{11}}\right) - 1 \,.
}$

\end{enumerate}

\end{enumerate}

The  bifurcations found numerically are in agreement with the qualitative findings of Proposition \ref{summary_prop}.  Also,  the EM diagrams for the $H_3$-type, and the $D_2H$ and $H_2D$-types molecules with chosen parameters fulfilling  the   $2 (b)$ case conditions  of Theorem (\ref{summary_prop}) are qualitatively in  good  agreement with the diagrams of \cite{Ko99} and \cite{Ko00}. 
Note that by  Remark \ref{non_stab_crit}  isosceles RE are never Lyapunov  stable.  




\begin{figure}[h!]
\center
  \includegraphics[scale=0.4]{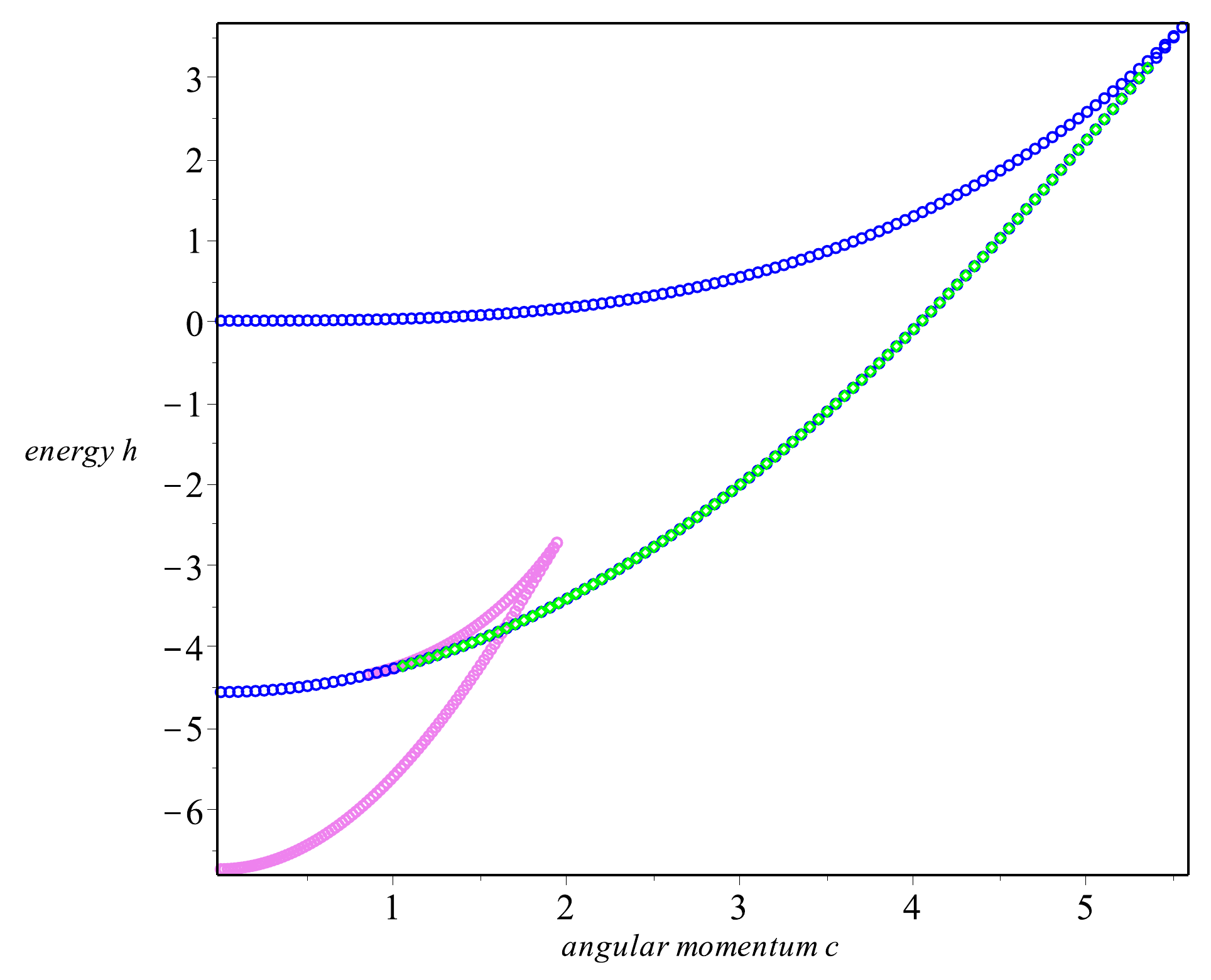}   
   \caption{The  EM diagram for $H_3$-type molecules, in agreement with the case $2 (b)$ of Proposition \ref{summary_prop} (see Remark \ref{remark_3_ident_atoms}). The unstable and stable planar linear RE are depicted in blue and green, respectively. The isosceles RE, all unstable, are depicted in violet (color online). We used  $m_A=m_B=1$, $a_{11}=a_{22}=3, \,\,b_{11}=b_{22}=1$.}
   \label{REM_H3}
   \end{figure}

   \begin{figure}[h!]
\centering
      \subfigure [The  EM diagram  for $D_2H$-type molecules in agreement with the case $2 (a)$ of Proposition \ref{summary_prop}. We used  $m_A=1,$ $m_B=0.5$, $a_{11}=6,$ $a_{22}=400, \,\,b_{11}=5,$ $b_{22}=3$.]
       {
            \includegraphics[scale=0.322]{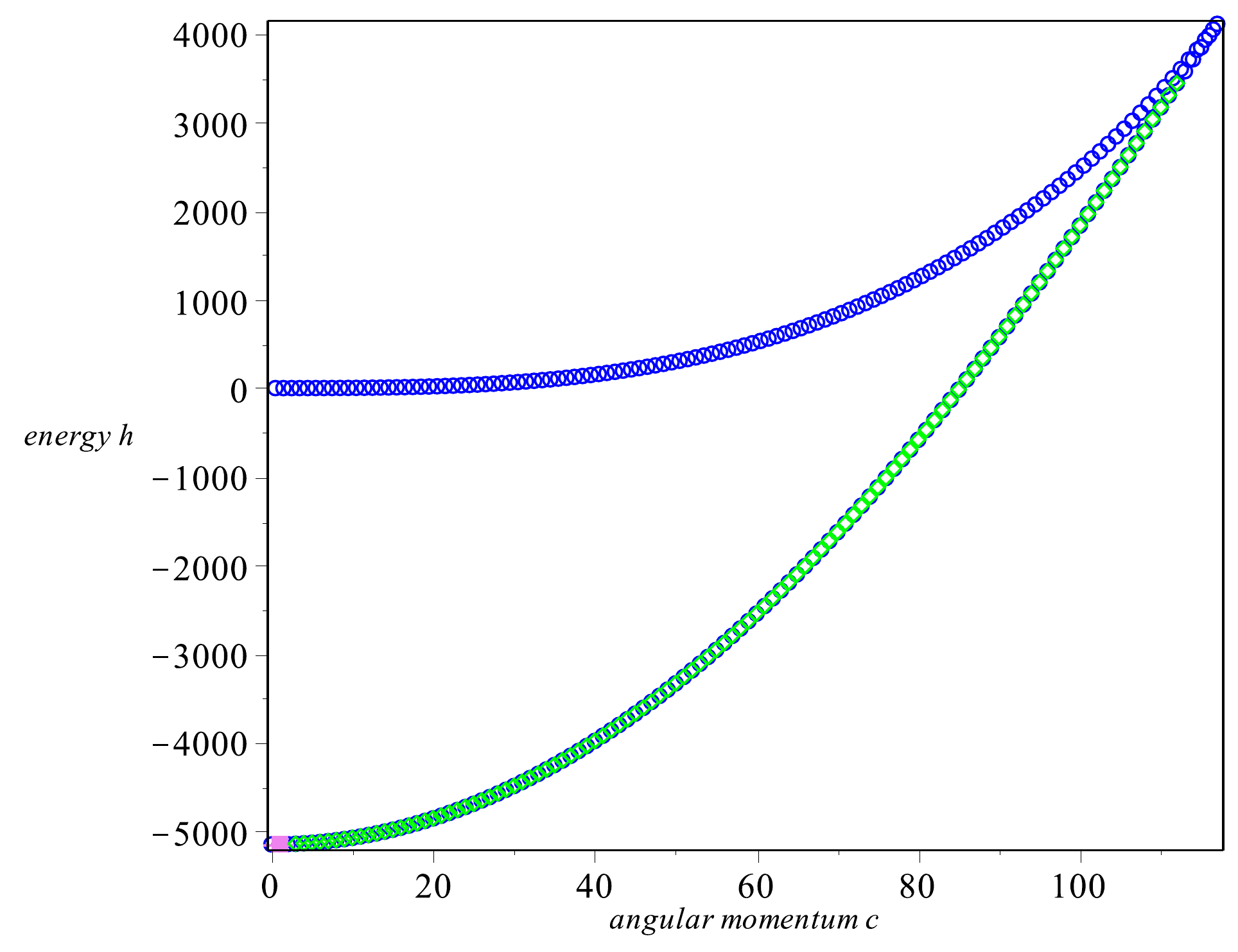}
         } 
         \hspace{2cm} 
         \subfigure [The  EM for $D_2H$-type molecules in agreement with the case $2 (b)$ of Proposition \ref{summary_prop}. We used   $m_A=1$, $m_B=0.5$  and $a_{11}= 4 ,\,\,a_{22}=2,      \,\,b_{11}=3,    \,\,b_{22}= 1$.]
       {
            \includegraphics[scale=0.3]{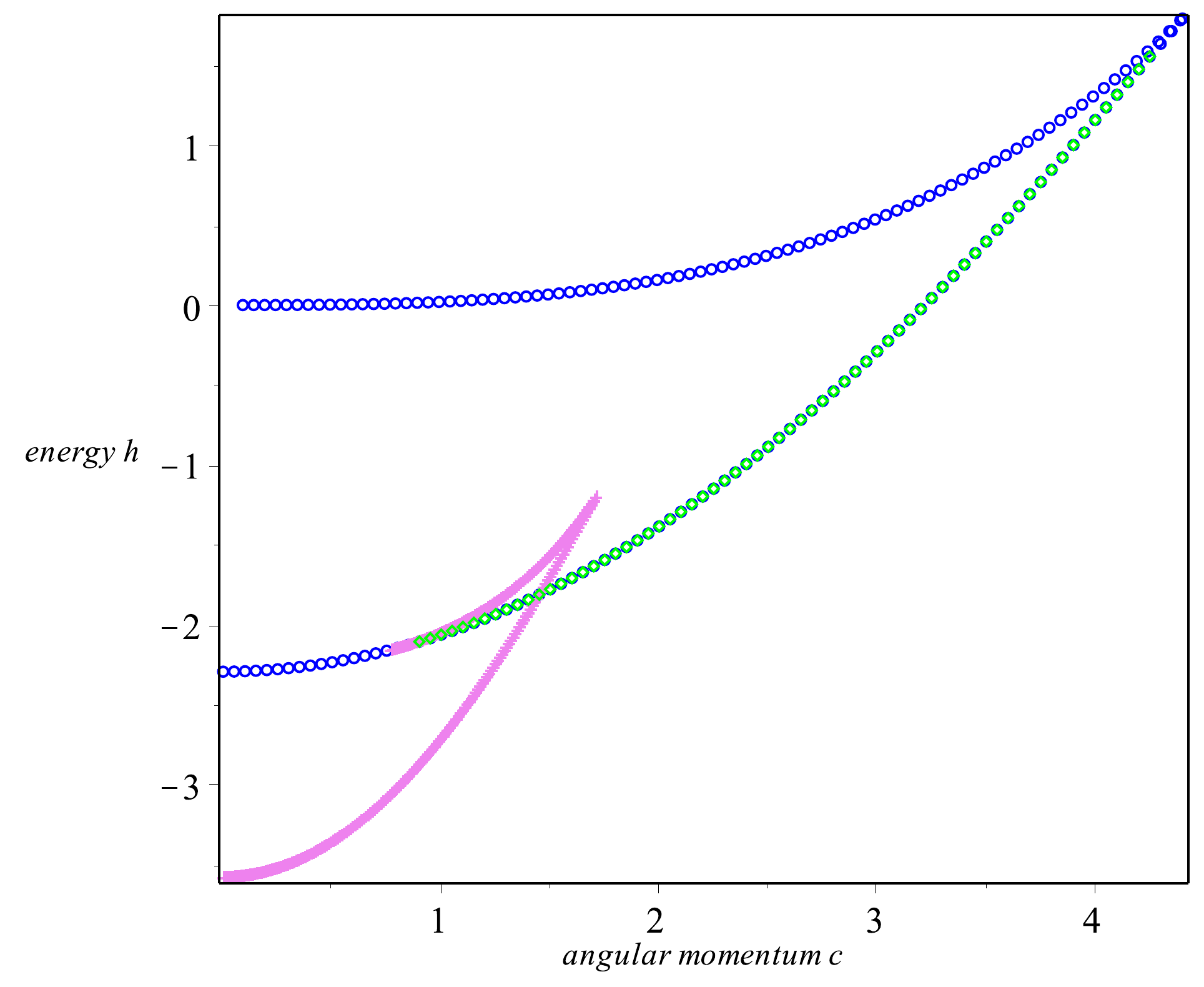}
         } 
         \caption{The  EM diagram for $D_2H$-type molecules. Unstable and stable planar linear RE are depicted in blue and green, respectively. The isosceles RE, all unstable, are depicted in violet (color online).
         } 
         \label{REM_D2H}
\end{figure}

   \begin{figure}[h!]
\centering
      \subfigure [The  EM diagram for $H_2D$-type molecules in agreement with the case $2 (a)$ of Proposition \ref{summary_prop}. We used  $m_A=0.5,$ $m_B=1$, $a_{11}=3,$ $a_{22}=200, \,\,b_{11}=1,$ $b_{22}=2$.]
       {
            \includegraphics[scale=0.32]{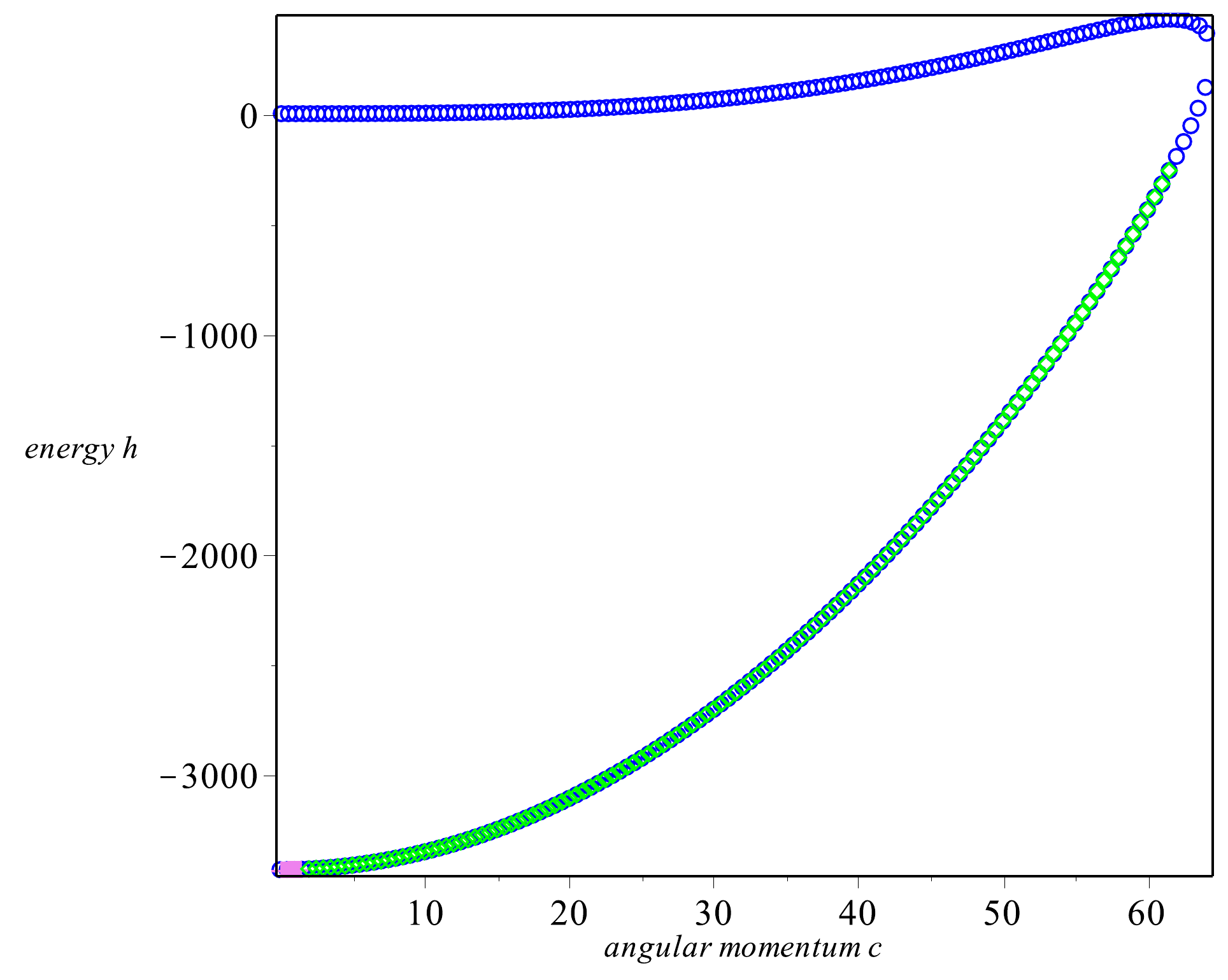}
         } 
         \hspace{1cm} 
         \subfigure [The  EM diagram for $H_2D$-type molecules in agreement with the case $2 (b)$ of Proposition \ref{summary_prop}. We used   $m_A=0.5$, $m_B=1$  and $a_{11}= 2 ,\,\,a_{22}=1,      \,\,b_{11}=1,    \,\,b_{22}= 3$.) ]
       {
            \includegraphics[scale=0.323]{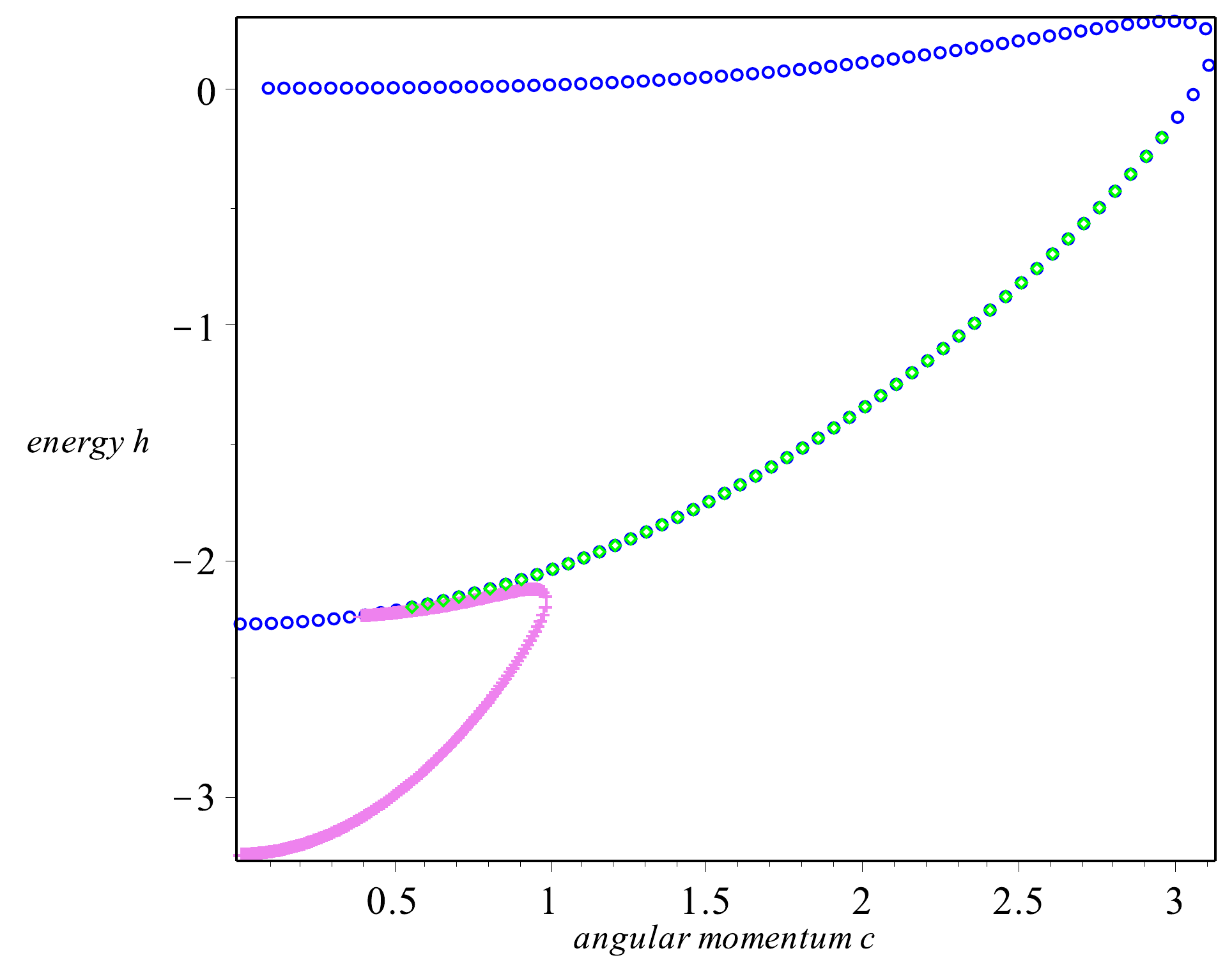}
         } 
         \caption{The  EM diagram for $H_2D$-type molecules. Unstable and stable planar linear RE are depicted in blue and green, respectively (color online). The isosceles RE, all unstable, are depicted in violet (color online).   
         } 
         \label{REM_H2D}
\end{figure}

\section{Conclusions}
In this paper we study isosceles triatomic ``A-B-A" molecules rotating about their symmetry axis. We show that the existence and bifurcations of the rotating steady state solutions may be determined qualitatively assuming the knowledge of shape of the diatomic A-A and A-B potentials  and three associated measurable parameters.
Consequently, our analysis shows that the RE bifurcations of an A-B-A triatomic molecule are not specified solely by the shape of the diatomic binary potentials A-A and A-B, but  additional parameters  are required.

We further present two geometric mechanics methods for determining stability. For isosceles RE we employ the REM method of Marsden and co-workers \cite{Ma92}, whereas for linear RE we implement the less-known Symplectic Slice method of Roberts \&al. \cite{RSS06}. We verify our qualitative results and calculate RE stability on some numerical examples using Lennard-Jones -type potentials.

\section{Acknowledgements}
CS  work was supported by a NSERC Discovery grant. DM completed part of this work during a MSc program at Wilfrid Laurier University. The authors thank the anonymous referee for his/her remarks and suggestions.


\end{document}